\documentclass[12pt,reqno]{amsart}

\usepackage{amssymb}

\newcommand{\lcm}{\operatorname{lcm}}
\newcommand{\ord}{\operatorname{ord}}
\newcommand{\fix}{\operatorname{\mathcal F}}
\newcommand{\mertens}{\operatorname{\mathcal M}}

\newcommand{\orbit}{\operatorname{\mathcal O}}
\def\divides{{\mathchoice{\mathrel{\bigm|}}{\mathrel{\bigm|}}{\mathrel{|}}{\mathrel{|}}}}
\def\notdivides{\mathrel{\kern-3pt\not\!\kern3.5pt\bigm|}}
\def\smallnotdivides{\mathrel{\kern-2pt\not\!\kern3.5pt\vert}}
\def\Divides{\divides\!\divides}
\def\bigo{\operatorname{O}}    
\def\littleo{\operatorname{o}} 
\newcommand{\dee}{\operatorname{d}}
\def\cP{{\mathcal P}}
\def\bt{{\boldsymbol t}}
\def\bp{{\boldsymbol p}}
\def\br{{\boldsymbol r}}
\def\bs{{\boldsymbol s}}

\def\bO{{\boldsymbol 0}}
\def\fp{v}

\renewcommand{\le}{\leqslant}
\renewcommand{\ge}{\geqslant}

\renewcommand{\mid}{\colon}
\newtheorem{theorem}{Theorem}[section]
\newtheorem{proposition}[theorem]{Proposition}
\newtheorem{lemma}[theorem]{Lemma}
\newtheorem{corollary}[theorem]{Corollary}

\theoremstyle{definition}

\newtheorem{example}[theorem]{Example}

\theoremstyle{remark}
\newtheorem{remark}[theorem]{Remark}

\begin{document}

\bibliographystyle{plain}

\title{Orbit-counting
in non-hyperbolic dynamical systems}

\author{G. Everest}
\author{R. Miles}
\author{S. Stevens}
\author{T. Ward}

\dedicatory{Draft \today}
\address{School of Mathematics, University of East
  Anglia, Norwich NR4 7TJ, United Kingdom}
\email{t.ward@uea.ac.uk}

\thanks{This research was supported by E.P.S.R.C. grant EP/C015754/1.}

\renewcommand{\subjclassname}{\textup{2000} Mathematics Subject Classification}
\begin{abstract}
There are well-known analogs of the prime number theorem and
Mertens' theorem for dynamical systems with hyperbolic behaviour.
Here we consider the same question for the simplest non-hyperbolic
algebraic systems. The asymptotic behaviour of the orbit-counting
function is governed by a rotation on an associated compact group,
and in simple examples we exhibit uncountably many different
asymptotic growth rates for the orbit-counting function. Mertens'
Theorem also holds in this setting, with an explicit rational
leading coefficient obtained from arithmetic properties of the
non-hyperbolic eigendirections.
\end{abstract}

\maketitle

\section{Introduction}

A closed orbit~$\tau$ of length~$\vert\tau\vert=n$ for a
continuous
map~$T:X\to X$ is a set of the form~$
\{x,T(x),T^2(x),\dots,T^n(x)=x\}$ with cardinality~$n$. A dynamical
analog of the prime number theorem concerns the asymptotic behaviour
of expressions like
\begin{equation}\label{equation:orbitvcounter}
\pi_T(N)=\left\vert\{\tau\mid\vert\tau\vert\le N\}\right\vert,
\end{equation}
and a dynamical analog of Mertens' Theorem concerns asymptotic
estimates for expressions like
\begin{equation}\label{equation:mertensequation}
\mertens_T(N)=\sum_{\vert\tau\vert\le
N}\frac{1}{e^{h(T)\vert\tau\vert}}
\end{equation}
where~$h(T)$ denotes the topological entropy of the map. Results
about the asymptotic behaviour of both expressions under the
assumption that~$X$ has a metric structure with respect to which~$T$
is hyperbolic may be found in the works of Parry~\cite{MR85c:58089},
Parry and Pollicott~\cite{MR85i:58105}, Sharp~\cite{MR93a:58142} and
others. An orbit-counting result on the asymptotic behavior
of~\eqref{equation:orbitvcounter} for quasi-hyperbolic toral automorphisms has been
found by Waddington~\cite{MR92k:58219}, and an
analog of Sharp's dynamical Mertens' Theorem for quasi-hyperbolic toral automorphisms has been
found by Noorani~\cite{MR1787646}. Both the current state of
these kinds of results and the seminal early work on geodesic flows
is described in the book of Margulis~\cite{MR2035655} which also has
a survey by Sharp on periodic orbits of hyperbolic flows.

One of the tools used in studying orbit-growth properties of
hyperbolic maps is the dynamical zeta function. This may be viewed
as a generalization of the Weil zeta function, which corresponds to
the dynamical zeta function of the action of the Frobenius map on
the extension of an algebraic variety over a finite field to the
field's algebraic closure. Writing
\begin{equation*}
\fix_T(n)=\vert\{x\in X\mid T^nx=x\}\vert
\end{equation*}
for the number of points fixed by~$T^n$, the dynamical zeta function
is defined by the formal expression
\begin{equation}\label{equation:dynamicalzetafunction}
\zeta_{T}(z)=\exp\sum_{n=1}^{\infty}\frac{z^n}{n}\fix_T(n)
\end{equation}
which has a formal expansion as an Euler product,
\begin{equation}\label{equation:eulerexpansion}
\zeta_T(z)=\prod_{\tau}\left(1-z^{\vert\tau\vert}\right)^{-1},
\end{equation}
where the product is taken over all orbits of~$T$. Just as the
classical Euler product relates analytic properties of the Riemann
zeta function to asymptotic counting properties of the prime
numbers, the Euler expansion~\eqref{equation:eulerexpansion} relates
analytic properties of the dynamical zeta function to orbit-counting
asymptotics. In the hyperbolic case, the zeta
function~\eqref{equation:dynamicalzetafunction} has radius of
convergence~$e^{-h(T)}$ and, crucially, has a meromorphic extension
to a strictly larger radius.

Our purpose here is on the one hand to study a very special class of
maps of arithmetic origin, while on the other relaxing the
hyperbolicity or quasi-hyperbolicity assumption. In this setting,
the simplest non-trivial example is the map~$\phi:X\to X$ dual to
the map~$r\mapsto 2r$ on~$\mathbb Z[\frac{1}{3}]$. This map is an
isometric extension of the circle-doubling map~$\psi(t)=2t\pmod{1}$
on the additive circle~$\mathbb T$ by a cocycle taking values in
the~$3$-adic integers~$\mathbb Z_3$; it is non-expansive and has
topological entropy~$\log2$. The dynamical zeta-function associated
to the map~$\phi$ is shown to have a natural boundary by Everest,
Stangoe and Ward~\cite{esw}, making it impossible to find a
meromorphic extension beyond the radius of convergence. The radius
of convergence is~$e^{-h(\phi)}=\frac{1}{2}$ since easy estimates
show that
\[
\frac{1}{n}\log\fix_{\phi}(n)\rightarrow\log2\mbox{ as }n\to\infty.
\]
The
bounds
\begin{equation}\label{Youarethetigerburningbright}
\frac{1}{3}\le\liminf_{N\to\infty}\frac{N\pi_{\phi}(N)}{2^{N+1}}
\le\limsup_{N\to\infty}\frac{N\pi_{\phi}(N)}{2^{N+1}}\leqslant1
\end{equation}
were found in~\cite{esw}. A problem left open there is to describe
the asymptotics exactly, and in particular to show
that~$\frac{N\pi_{\phi}(N)}{2^{N+1}}$ does not converge
as~$N\to\infty$.

A similar result is found for the dynamical analog of Mertens'
Theorem. Write
\[
\orbit_{T}(n)=\left\vert\{\tau\mid \tau \mbox{ is a closed orbit
of~$T$ of length }\vert\tau\vert=n\}\right\vert
\]
for the number of orbits of length~$n$ under~$T$. Then
\begin{equation}\label{Deepintheforestofmynight}
\frac{1}{2}\log{N}+\bigo(1)\leqslant \sum_{n\leqslant
N}\frac{\orbit_{\phi}(n)}{2^n} \leqslant\log N+\bigo(1)
\end{equation}
is shown in~\cite{esw}.

A consequence of the results in this paper is a better explanation
of the sequences along which the expressions
in~\eqref{Youarethetigerburningbright} converge, and a proof that
there is a single asymptotic in~\eqref{Deepintheforestofmynight}.
The map considered in~\cite{esw} is a special case of a more general
construction of~$S$-integer maps described in~\cite{MR99b:11089}.
These are parameterized by an~$\mathbb A$-field~$\mathbb K$ (for
example,~$\mathbb Q$ or~$\mathbb F_q(t)$), a subset~$S$ of the set
of places of~$\mathbb K$, and an element~$\xi\in\mathbb K^*$ of
infinite order (see the start of
Section~\ref{section:proofofmoregeneralorbittheorem} for the
construction; the assumption that~$\xi$ has infinite multiplicative
order is equivalent to ergodicity for the resulting map). For the
map~$\phi$ above, these parameters are chosen with~$\mathbb
K=\mathbb Q$,~$S=\{3\}\subset\{2,3,5,7,11,\dots\}$ and~$\xi=2$. If
the~$\mathbb A$-field~$\mathbb K$ has characteristic zero, then the
resulting map is an endomorphism of a solenoid.

The essential starting point is to note
from~\cite{MR99b:11089}
that if~$T:X\to X$ is
an~$S$-integer map with~$S$ finite and~$X$ connected, then
\[
    \frac{1}{n}\log\fix_T(n)\longrightarrow h(T)>0,
\]
so the dynamical zeta function has radius of
convergence~$e^{-h(T)}$. This suggests that the natural function to
compare~$\pi_T(N)$ with is~$\frac{e^{h(T)(N+1)}}{N}$, so
define
\[
\Pi_T(N)=\frac{N\pi_{T}(N)}{e^{h(T)(N+1)}}.
\]

\begin{theorem}\label{theorem:moregeneralorbitcounting}
Let~$T:X\to X$ be an~$S$-integer map with~$X$ connected and~$S$
finite.
Then~$\left(\Pi_T(N)\right)$ is a bounded sequence,
and
\[
\liminf_{N\to\infty}\Pi_T(N)>0.
\]
Moreover, there is an
associated pair~$(X^*,a_T)$, where~$X^*$ is a compact group
and~$a_T\in X^*$, with the property that if~$a_T^{N_j}$ converges
in~$X^*$ as~$j\to\infty$,
then~$\Pi_T(N_j)$ converges in~$\mathbb
R$ as~$j\to\infty$.
\end{theorem}

Thus the pair~$(X^*,a_T)$ detects limit points in the orbit-counting
problem. In the hyperbolic case, the group~$X^*$ is trivial,
reflecting the fact that~$\left(\Pi_T(N)\right)_{N\ge1}$ itself
converges.

\begin{example}\label{example:quasihyperbolic}
The most familiar examples of non-hyperbolic automorphisms are the
quasi-hyperbolic toral automorphisms (see Lind~\cite{MR84g:28017}
for a detailed account of their dynamical properties.)
Let~$k=\mathbb Q(\xi)$ where~$
\xi=-(1+\root\of{2})-\root\of{2\sqrt{2}+2},$ and~$S=\emptyset$. Then
the corresponding map~$T$ is the quasi-hyperbolic automorphism of
the~$4$-torus defined by the matrix
\[
\left[
\begin{matrix}0&1&0&0\\0&0&1&0\\0&0&0&1\\-1&-4&2&-4\end{matrix}
\right].
\]
There is a pair of eigenvalues~$\lambda,\overline{\lambda}$
with~$\vert\lambda\vert=1$. The corresponding system~$(X^*,a_T)$ is
the rotation~$z\mapsto\lambda z$ on~$\mathbb S^1$, and any
sequence~$(N_j)$ for which~$\left(\lambda^{N_j}\right)$ converges
has the property that~$\left(\Pi_T(N_j)\right)$
converges as~$j\to\infty$. This recovers in part a result of
Waddington~\cite{MR92k:58219}, who explicitly
identifies~$\Pi_T(N)$ as an almost-periodic
function of~$N$.
\end{example}

In some cases the correspondence between convergent subsequences
seen in the detector group~$X^*$ and the orbit-counting problem is
exact. For simplicity we state this for the case~$\mathbb K=\mathbb
Q$,~$\xi=2$,~$S=\{3\}$; the same method gives a similar conclusion
whenever~$\mathbb K=\mathbb Q$ and~$\vert S\vert=1$. The full extent
of the phenomena (and, in particular, of the appearance of
uncountably many limit points) is not clear.

\begin{theorem}\label{theorem:23example}
For the map~$\phi$ dual to the map~$x\mapsto 2x$ on~$\mathbb
Z[\frac{1}{3}]$, the sequence~$\left(\Pi_{\phi}(N_j)\right)$
converges as~$j\to\infty$ if and only if the
sequence~$\left(2^{N_j}\right)$ converges in the group~$\mathbb
Z_3$. In particular, the sequence~$\left(\Pi_{\phi}(N)\right)$ has
uncountably many limit points. Moreover, the upper and lower limits
are both transcendental.
\end{theorem}

The dynamical analog of Mertens Theorem concerns the
expression~\eqref{equation:mertensequation}. In the simplest case
(an endomorphism of a~$1$-dimensional solenoid) precise results are
readily found, with a rational coefficient of the leading term.

\begin{theorem}\label{theorem:mertensinsimplecases23}
For an~$S$-integer map~$T$ corresponding to~$\mathbb K=\mathbb Q$
and~$S$ finite, there are constants~$k_T\in\mathbb Q$ and~$C_T$ such
that
\[
\mertens_{T}(N)=k_T\log N+C_T+\bigo\left(1/N\right).
\]
\end{theorem}

\begin{example}\label{example:explicitconstants}
Let~$\xi=2$ in Theorem~\ref{theorem:mertensinsimplecases23}, so the
map~$T$ is the map dual to~$x\mapsto 2x$ on the ring~$
R_S=\{\textstyle\frac{p}{q}\in\mathbb Q\mid\mbox{ primes dividing
}q\mbox{ lie in }S\}.$ The constant~$k_T$ for various simple sets~$S$
is given in Table~\ref{table}.
\begin{table}[ht!]
\begin{center}
\caption{\label{table}Leading coefficients in Mertens' Theorem}
\begin{tabular}{c|c}
$S$&value of $k_T$\\
\hline
$\emptyset$&$1$\\
$\{3\}$&$\frac{5}{8}\vphantom{\displaystyle\sum}$\\
$\{3,5\}$&$\frac{55}{96}\vphantom{\displaystyle\sum}$\\
$\{3,7\}$&$\frac{269}{576}\vphantom{\displaystyle\sum}$\\
co-finite&$0$\\
\end{tabular}
\end{center}
\end{table}
\end{example}

In the general case there is less control of the error term
(the error term in the dynamical Mertens' Theorem of Sharp~\cite{MR93a:58142}
for the hyperbolic setting
is improved to~$\littleo(1/N)$ by Pollicott~\cite{MR1045147}).

\begin{theorem}\label{theorem:mertensforSfinite}
Let~$T:X\to X$ be an~$S$-integer map with~$X$ connected and with~$S$
finite. Then there are constants~$k_T\in\mathbb Q$,~$C_T$
and~$\delta>0$ with
\begin{equation}\label{equation:threetermestimate}
\mertens_T(N)=k_T\log N+C_T+\bigo(N^{-\delta}).
\end{equation}
\end{theorem}

At the other extreme, the class of~$S$-integer systems with~$\vert
S\vert$ infinite provides a range of subtle behaviors that cannot
readily be treated in this way. Possibilities include~$\fix(n)$
growing much slower than exponentially; the `generic' behavior
for~$S$ chosen randomly is discussed in~\cite{MR1458718}
and~\cite{MR99k:58152}. Some results on systems with~$S$ co-finite
may be found in the thesis of Stangoe~\cite{stangoethesis}.

\begin{example}
Let~$T$ be an~$S$-integer map dual to~$x\mapsto \xi x$ with~$\mathbb
K=\mathbb Q$ and~$S$ co-finite. For any finite place~$w\in S$ there
are constants~$A,B>0$ with~$\vert\xi^n-1\vert_w>A/n^B$, so by the
product formula there is a constant~$C>0$ with~$\fix_{T}(n)\le n^C$.
It follows that~$\mertens_{T}(N)$ is bounded for all~$N$.
\end{example}

Allowing the compact group~$X$ to be infinite-dimensional is
problematical for a different reason: the following example may be
found in~\cite[Th.~8.1]{stangoethesis}.

\begin{example}
For any sequence~$a_1,a_2,\dots$ there is an automorphism~$T$ of a
compact connected group with
\[
a_n\le\fix_T(n)<\infty\mbox{ for all }n\ge1.
\]
To see this, define a sequence of maps~$T_1,T_2,\dots$ as follows.
Let~$T_{1}$ be the map dual to~$x\mapsto3x$ on~$\mathbb{Z}$.
Let~$T_{2}$ be the map dual to~$x\mapsto2x$ on~$\mathbb{Z}$.
Let~$T_{3}$ be the map dual to~$x\mapsto3x$
on~$\mathbb{Z}[\frac{1}{2}]$. By Zsigmondy's Theorem,
\[
\{p\mid p\divides 3^{n}-1\mbox{ for some }n\leq k\}\subsetneq \{p
\mid p\divides 3^{n}-1\mbox{ for some }n\leq k+1\}
\]
unless~$k=1$. This allows the sequence of maps to be continued:
Let~$T_{4}$ be the map dual to~$x\mapsto3x$
on~$\mathbb{Z}[\frac{1}{2},\frac{1}{13}]$ and, similarly~$T_{k}$
will be the map dual to~$x\mapsto3x$
on~$\mathbb{Z}[\frac{1}{s_1},\dots,\frac{1}{s_t}]$, where
\[
\{s_1,\dots,s_t\}=\{p\mid p\mbox{ is a prime with }p\divides
3^n-1\mbox{ for some }n< k\}.
\]
Using the periodic point
formula~\eqref{equation:periodicpointformulaingeneral}
from~\cite{MR99b:11089}, the choice of primes ensures
that~$\fix_{T_k}(j)=1$ for~$j<k$ and~$\fix_{T_k}(k)>1$. Finally
define the map~$T$ to be the infinite product
\[
T=\underbrace{\left(T_{1}\times T_{1}\times\cdots\times
T_{1}\right)}_{\text{so that }\infty>\fix_{T}(1)>a_{1}}\times
\underbrace{\left(T_{2}\times T_{2}\times\cdots\times
T_{2}\right)}_{\text{so that }\infty >\fix_{T}(2) > a_{2}}\times
\cdots.
\]
For any~$k\ge1$, all but finitely many terms in the product
giving~$\fix_T(n)$ are~$1$, so the product is finite and
exceeds~$a_n$.
\end{example}

The paper is organized as follows. Theorem~\ref{theorem:23example}
and Theorem~\ref{theorem:mertensinsimplecases23} for the same
map~$\phi$ dual to~$x\mapsto2x$ on~$\mathbb Z[\frac{1}{3}]$ are proved in
Section~\ref{section:boththeoremsfor23example}; this example
illustrates some of the issues that arise in the more general
setting while avoiding the Diophantine subtleties.
Theorem~\ref{theorem:moregeneralorbitcounting} is proved in
Section~\ref{section:proofofmoregeneralorbittheorem}.
Theorem~\ref{theorem:mertensforSfinite} without an error term is
proved in Section~\ref{section:mertenswithkequalsq}; this result may
be found using soft methods.
Theorem~\ref{theorem:mertensinsimplecases23} is proved in
Section~\ref{section:mertenswithkequalsq}, with the essential
combinatorial step generalized to allow other fields. Finally,
Section~\ref{section:allowinginfiniteplaces} assembles the
additional Diophantine ingredients for
Theorem~\ref{theorem:mertensforSfinite}.

\section{Proof of Theorems~\ref{theorem:23example} and~\ref{theorem:mertensinsimplecases23}
in a special case} \label{section:boththeoremsfor23example}

The specific map~$\phi$ dual to~$x\mapsto 2x$ on~$\mathbb
Z[\frac{1}{3}]$ already reveals some of the essential features of
these systems. In addition, the relatively simple nature of the map
allows very precise results. This section contains a self-contained
proof of Theorem~\ref{theorem:23example} which may be read on its
own or used to motivate some of the arguments in
Section~\ref{section:proofofmoregeneralorbittheorem}. It also
contains a self-contained proof of
Theorem~\ref{theorem:mertensinsimplecases23}
for the case~$S=\{3\}$ and~$\xi=2$.

By~\cite[Lem.~5.2]{MR99b:11089}, the number of points fixed
by~$\phi^n$ is
\[
\fix_{\phi}(n)=(2^n-1)\vert 2^n-1\vert_3,
\]
so the number of orbits of length~$n$ is given by
\[
\orbit_{\phi}(n)={\frac{1}{n}}\sum_{d\vert n}\mu(
{\textstyle\frac{n}{d}})(2^d-1)|2^d-1|_3
\]
by M{\"o}bius inversion, and hence
\begin{equation}\label{Downinthedarknessofmydreams}
\pi_{\phi}(N)=\sum_{n\le N}{\frac{1}{n}}\sum_{d\vert n}\mu(
\textstyle\frac{n}{d})(2^d-1)|2^d-1|_3.
\end{equation}
We begin by replacing~\eqref{Downinthedarknessofmydreams} with a
more manageable expression. Let
\begin{equation}\label{equation:GNfirstdefined}
G(N)=\sum_{n\le N}{\frac{1}{n}} 2^n\vert 2^n-1\vert_3.
\end{equation}
Then
\begin{eqnarray*}
\left\vert\pi_{\phi}(N)-G(N)\right\vert&\le& \sum_{n\le
N}{\frac{1}{n}}\left( \sum_{d\divides n}
\underbrace{\vert2^d-1\vert_3}_{\le1}+
\sum_{d\divides n,d<n}2^d\vert 2^d-1\vert_3\right)\\
&\le& \sum_{n\le N}{\frac{1}{n}} \left( n+\sum_{d\le\lfloor
n/2\rfloor}2^d\right)= \bigo\left(2^{N/2}\right),
\end{eqnarray*}
so for the purposes of the asymptotic sought we can use~$G(N)$ in
place of~$\pi_{\phi}(N)$.

We next give a simple proof of the orbit-counting asymptotic for the
circle-multiplication by~$a\ge 2$, that is for the
map~$\psi_a(x)=ax\pmod{1}$; for this map~$\fix_{\psi_a}(n)=a^n-1$.
Results like these are special cases of the more general picture in
the work of Parry and Pollicott~\cite{MR85i:58105}. We give an
elementary proof here because the argument used presages the
estimates needed later.

\begin{lemma}\label{YouwerethedreamIdseecometrue}
$\pi_{\psi_a}(N)\sim\frac{a^{N+1}}{N(a-1)}.$
\end{lemma}

\begin{proof}
By M{\"o}bius inversion
\[
\pi_{\psi_a}(N)=\sum_{n\le N}{\frac{1}{n}}\sum_{d\divides n}
\mu({\textstyle\frac{n}{d}})(a^d-1)= \sum_{n\le
N}{\frac{1}{n}}\sum_{d\divides n} \mu({\textstyle\frac{n}{d}})a^d-1.
\]
Subtracting the dominant terms,
\begin{eqnarray}
\left\vert \pi_{\psi_a}(N)-\sum_{n\le N}{\frac{1}{n}}a^n
\right\vert&=& 1+\left\vert \sum_{n\le N}{\frac{1}{n}}
\sum_{d\divides n,d<n}\mu({\frac{n}{d}})a^d
\right\vert\nonumber\\
&=& \bigo\left( \sum_{n\le N}\sum_{d\le\lfloor n/2\rfloor}a^d
\right)\nonumber\\
&=&\bigo\left(\sum_{n\le N}a^{n/2}\right)=
\bigo\left(a^{N/2}\right).\label{IwastheonewholookedsohardIcouldnotsee}
\end{eqnarray}
To estimate the dominant terms, let~$K(N)=\lfloor N^{1/4}\rfloor$.
Then
\begin{eqnarray*}
\left\vert \sum_{n\le N}{\frac{1}{n}}a^n-\sum_{N-K(N)\le n\le
N}{\frac{1}{n}}a^n\right\vert&\le&
\sum_{n\le N-K(N)}a^n\\
&=&\bigo\left( a^{N-K(N)}\right).
\end{eqnarray*}
Finally,
\begin{eqnarray*}
\sum_{N-K(N)\le n\le N}{\frac{1}{n}}a^n &=&
\frac{a^N}{N}\sum_{r=0}^{K(N)}a^{-r}\left(1-
{\frac{r}{N}}\right)^{-1}\\
&=&\frac{a^N}{N}\left[\vphantom{\sum}\right. \frac
a{a-1}-\bigo(a^{-K(N)})+\bigo(\sum_{r=0}^{K(N)}r/N)
\left.\vphantom{\sum}\right]\\
&=& \frac{a^{N+1}}{N(a-1)}+\bigo\left(
\frac{a^N}{N^2}\sum_{r=0}^{K(N)}r\right)\\
&=& \frac{a^{N+1}}{N(a-1)}+\bigo\left(\frac{a^N}{N^{3/2}}\right).
\end{eqnarray*}
Together with~\eqref{IwastheonewholookedsohardIcouldnotsee}, this
proves the lemma.
\end{proof}

Returning to the main problem, write
\[
I(N)=\sum_{n\le N,2\divides n} {\frac{1}{n}}2^n\vert2^n-1\vert_3
\]
and
\[
J(N)=\sum_{n\le N,2\smallnotdivides n}
{\frac{1}{n}}2^n\vert2^n-1\vert_3,
\]
so~$G(N)=I(N)+J(N)$. Splitting into odd and even terms further
simplifies the expressions since an easy calculation shows that
\begin{equation}\label{KatieMelua}
\vert2^n-1\vert_3= \left\{
\begin{array}{cl}
\frac{1}{3}\vert n\vert_3&\mbox{if $n$ is even;}\\
1&\mbox{if $n$ is odd,}
\end{array}
\right.
\end{equation}
so
\[
J(N)=\sum_{n\le N,2\smallnotdivides n} {\frac{1}{n}}2^n.
\]

\begin{lemma}
$J(N)\sim\frac{1}{3}\cdot\frac{2^{N+1}}{N}.$
\end{lemma}

\begin{proof}
Lemma~\ref{YouwerethedreamIdseecometrue} applied to the
maps~$\psi_2$ and~$\psi_4$ shows that
\[
\sum_{n\le N} {\frac{1}{n}}2^n\sim\frac{2^{N+1}}{N}\mbox{ and
}\sum_{k\le K}\frac{1}{k}4^k\sim\frac{4^{K+1}}{3K}.
\]
Hence
\begin{eqnarray*}
J(N)&=&\sum_{n\le N} {\frac{1}{n}}2^n- \sum_{n\le
N,2\divides n} {\frac{1}{n}}2^n\\
&=& \sum_{n\le N} {\frac{1}{n}}2^n-\sum_{k\le N/2}
{\frac{1}{2k}}4^k\\
&\sim& \frac{2^{N+1}}{N}-\frac{2}{3}\cdot \frac{2^{N+1}}{N}\\
&=&\frac{1}{3}\cdot\frac{2^{N+1}}{N}.
\end{eqnarray*}
\end{proof}

We are therefore left with the expression
\begin{eqnarray*}
I(N)=\sum_{n\le N,2\divides n} {\frac{1}{n}}2^n\vert2^n-1\vert_3&=&
\frac{1}{3}\sum_{n\le N,2\divides n}
{\frac{1}{n}}2^n\vert n\vert_3\\
&=& \frac{1}{6}\sum_{k\le{N}/{2}} {\frac{1}{k}}2^{2k}\vert k\vert_3.
\end{eqnarray*}
Define
\[
L(M)=\sum_{n\le M}{\frac{1}{n}}4^n\vert n\vert_3
\]
and
\[
a_M=\frac{ML(M)}{4^M}.
\]
Again it is enough to look only at the large terms, since
\begin{eqnarray*}
\left\vert \sum_{M-K(M)\le n\le M}\frac{4^n}{n}\vert n\vert_3
-\sum_{n\le M}\frac{4^n}{n}\vert n\vert_3 \right\vert&\le&
\sum_{n\le K(M)}4^n=\bigo\left( 4^{K(M)} \right).
\end{eqnarray*}
Expanding from the last term gives
\begin{eqnarray*}
a_M&=&\frac{\vert M\vert _3}{1}+\frac{4^{-1}\vert M-1\vert
_3}{1-1/M}+
\frac{4^{-2}\vert M-2\vert _3}{1-2/M}+\cdots \\
&&\quad\quad\quad\quad\quad\quad\quad\quad\quad+\frac{4^{-K(M)}\vert
M-K(M)\vert_3}{1-K(M)/M}\\
&=& \frac{\vert M\vert _3}{1}+\frac{\vert M-1\vert _3}{4}+
\frac{\vert M-2\vert _3}{4^2}+ \dots +\frac{\vert M-K(M)\vert
_3}{4^{K(M)}}\\
&&\quad\quad\quad\quad\quad\quad\quad\quad\quad
+\bigo\left(\sum_{r=1}^{K(M)}r/M\right),
\end{eqnarray*}
and the error term is~$\bigo\left({M}^{-1/2}\right)$. Thus the limit
points mentioned in Theorem~\ref{theorem:23example} come from limit
points of the sequence~$(b_M)$ defined by
\begin{equation}\label{Ilivedlikeawildandlonelysoul}
b_M=\frac{\vert M\vert _3}{1}+\frac{\vert M-1\vert _3}{4}+
\frac{\vert M-2\vert _3}{4^2}+ \dots +\frac{\vert M-K(M)\vert
_3}{4^{K(M)}}.
\end{equation}
Clearly
\[
b_M\le 1+\frac{1}{4}+\frac{1}{4^2}+\dots=\frac{4}{3}
\]
and
\[
b_M\ge\frac{1}{4}
\]
because~$3\divides M$ implies that~$3\notdivides(M-1)$. These upper
and lower bounds imply upper and lower bounds of~$\frac{5}{9}$
and~$\frac{9}{24}$ respectively
in~\eqref{Youarethetigerburningbright}.

The shape of the expression~\eqref{Ilivedlikeawildandlonelysoul}
suggests that the lower limit will be seen along sequences highly
divisible by~$3$, and the upper limit along sequences not divisible
by~$3$, and this indeed turns out to be the case. To find limit
points, it is easier to work with the infinite sum rather
than~\eqref{Ilivedlikeawildandlonelysoul}, so notice first that if
\[
c_M=\sum_{j=0}^{\infty}\frac{\vert M-j\vert_3}{4^j}
\]
then~$\vert b_M-c_M\vert=\bigo\left(2^{K(M)-M}\right)$. Now
let~$\vert M_k\vert_3=3^{-k}$ so that (by the ultrametric
inequality)
\begin{eqnarray*}
c_{M_k}+t_k &=&\frac{1}{3^k}+\frac{\vert1\vert_3}{4}+
\frac{\vert2\vert_3}{4^2}+\frac{\vert3\vert_3}{4^3}+\cdots\\
&=& \frac{1}{3^k}+ \sum_{j=1}^{\infty}\frac{1}{4^j}-
\frac{2}{3}\sum_{j=1}^{\infty}\frac{1}{4^{3j}}
-\frac{2}{9}\sum_{j=1}^{\infty}\frac{1}{4^{9j}}
-\frac{2}{27}\sum_{j=1}^{\infty}\frac{1}{4^{27j}}
-\cdots\\
&=& \frac{1}{3^k}+\frac{1}{3}-
2\sum_{r=1}^{\infty}\frac{1}{3^r(4^{3^r}-1)}
\end{eqnarray*}
where
\[
t_k=\sum_{j=3^k}^{\infty}\frac{\vert j\vert_3-\vert
M-j\vert_3}{4^j}=\bigo(4^{-3^k}).
\]
Thus~$c_{M_k}$ converges
as~$k\to\infty$. Moreover, the limiting value is transcendental.

\begin{lemma}
The sum~$C=\displaystyle\sum_{r=1}^{\infty}\frac{1}{3^r(4^{3^r}-1)}$
is transcendental, and
\[
\liminf_{M\to\infty}c_M=\frac{1}{3}-2C.
\]
\end{lemma}

\begin{proof}
Let~$q_s=3^s(4^{3^s}-1)$. Then there is an integer~$p_s$ such that
\[
C_s=C-\frac{p_s}{q_s}=\sum_{r=s+1}^{\infty}
\frac{1}{3^r(4^{3^r}-1)}.
\]
Thus~$C_s=\bigo\left(3^{-s+1}4^{-3^{s+1}}\right)$, so
\[
0<\vert C-\frac{p_s}{q_s}\vert=\bigo\left( q_s^{-3}\right)
\]
showing that~$C$ is too well-approximable to be algebraic.

To see that this does give the lower limit, notice that
\[
c_{M_k}=\frac{1}{3^k}+\frac{1}{3}-2C-t_k.
\]
Any limit point along a sequence~$(M_k)$ with~$\ord_3(M_k)$ bounded
infinitely often is larger, and any limit point
with~$\ord_3(M_k)\to\infty$ must be this one.
\end{proof}

Essentially the same argument with~$M_k$ with~$\vert
M_k+1\vert_3=3^{-k}$ shows that
\[
\limsup_{M\to\infty}c_M=4\liminf_{M\to\infty}c_M,
\]
completing the proof of the first part of
Theorem~\ref{theorem:23example}.

We now turn our attention to the remaining part of
Theorem~\ref{theorem:23example}.

\begin{lemma}\label{Whatarealltheselongnames?}
Fix~$M,N\in\mathbb N$ with~$0<\varepsilon=\vert M-N\vert_3$. Then
\[
\frac{\varepsilon}{3\cdot4^{3/\varepsilon}} < \vert c_M-c_N\vert \le
\frac 43\varepsilon
\]
\end{lemma}

\begin{proof} The second inequality is straightforward: By the reverse
triangle inequality
\begin{equation} \label{Amoresensiblename}
\left\vert\vphantom{A^A}\right.\vert M-j\vert_3 - \vert N-j\vert_3
\left.\vphantom{A^A}\right\vert \le \vert M-N\vert_3 = \varepsilon
\end{equation}
for any~$j$, so that
\[
\vert c_M-c_N\vert \le  \sum_{j=0}^\infty
\frac{\left\vert\vphantom{A^A}\right.\vert M-j\vert_3 - \vert
N-j\vert_3 \left.\vphantom{A^A}\right\vert}{4^j} \le
\sum_{j=0}^\infty \frac{\varepsilon}{4^j} = \frac 43\varepsilon.
\]

For the first inequality a more careful analysis of where the series
in~$c_M$ and~$c_N$ differ is needed. Write~$\varepsilon=3^{-k}$,
with~$k\ge 0$. There exist unique integers~$0\le j_M,j_N< 3^{k+1}$
such that
\[
\vert M-j_M\vert_3 \le 3^{-(k+1)} \medspace\mbox{ and }\medspace
\vert N-j_N\vert_3 \le 3^{-(k+1)}.
\]
Since~$\vert M-N\vert_3=3^{-k}$ we have~$\vert
j_M-j_N\vert_3=3^{-k}$ also and we may assume that~$j_M<j_N$ without
loss of generality. By the ultrametric inequality,
\[
\vert M-j\vert_3 = \vert N-j\vert_3,\quad\mbox{ for }j<j_M,
\]
so the series in~$c_M$ and~$c_N$ differ first at the term~$j=j_M$.
Thus
\[
\vert M-j_M\vert_3\le 3^{-(k+1)}< \vert N-j_M\vert_3 =3^{-k}
\]
and so
\begin{eqnarray*}
\vert c_M-c_N\vert &\ge& \frac{\vert N-j_M\vert_3-\vert
M-j_M\vert_3}{4^{j_M}} -\negmedspace\negmedspace
\sum_{j=j_M+1}^\infty \! \frac{\left\vert\vphantom{A^A}\right. \vert
M-j\vert_3 - \vert N-j\vert_3
\left.\vphantom{A^A}\right\vert}{4^j} \\
&\ge& \frac{3^{-k}-3^{-(k+1)}}{4^{j_M}} - \frac{3^{-k}}{4^{j_M}}
\sum_{j=1}^\infty \frac 1{4^j} \\
&\ge& \frac{3^{-k}}{4^{j_M}}\left(1-\frac 13 -\frac 13\right) \ =\
\frac {3^{-k}}{3.4^{j_M}}\ >\
\frac{\varepsilon}{3.4^{3/\varepsilon}}
\end{eqnarray*}
by~\eqref{Amoresensiblename}.
\end{proof}

An immediate consequence of Lemma~\ref{Whatarealltheselongnames?} is
the following corollary, from which the remainder of
Theorem~\ref{theorem:23example} follows.

\begin{corollary}
Given any~$\alpha\in\mathbb Z_3$ and sequence of natural
numbers~$(M_k)$ converging to~$\alpha$ in~$\mathbb Z_3$,
define~$c_\alpha$ to be~$\lim_{M_k\to\infty} c_{M_k}$.
Then~$c_\alpha$ is well-defined (the limit exists and is independent
of the choice of approximating sequence). Moreover,
if~$\beta\in\mathbb Z_3$ and~$\varepsilon=\vert\alpha-\beta\vert_3$
then
\[
\frac{\varepsilon}{3.4^{3/\varepsilon}} \le \vert
c_\alpha-c_\beta\vert \le \frac 43\varepsilon.
\]
\end{corollary}

This completes the proof of Theorem~\ref{theorem:23example}.

Theorem~\ref{theorem:mertensinsimplecases23} for the map~$\phi$
concerns the sum
\begin{equation}\label{Mersum}
\mertens_{\phi}(N)=\sum_{n \le N}\frac{\orbit_{\phi}(n)}{2^n}
\end{equation}
where~$\orbit_{\phi}(n)$ is the number of orbits of length~$n$
under~${\phi}$, so
\[
\orbit_{\phi}(n)=\frac{1}{n}\sum_{d\vert
n}\mu({n}/{d})(2^d-1)|2^d-1|_3.
\]
Let
\begin{equation}\label{Youweretheonewhobroughtmehomedowntoearth}
F(N)=\sum_{n\le N}\frac{\vert2^n-1\vert_3}{n},
\end{equation}
and notice that
\begin{eqnarray*}
\mertens_{\phi}(N)-F(N)&=&\sum_{n\le N}\frac{1}{n}
\left(
\sum_{d\divides n}\mu({\textstyle{\frac{n}{d}}})\vert 2^d-1\vert_3\frac{2^d-1}{2^n}-\vert2^n-1\vert_3
\right)\\
&=&
\sum_{n\le N}\frac{1}{n}
\left(
\frac{\vert2^n-1\vert_3}{2^n}+\negmedspace
\sum_{d\divides n,d<n}
\mu({\textstyle{\frac{n}{d}}})
\vert2^d-1\vert_3\frac{2^d-1}{2^n}
\right)\\
&=&\sum_{n\le N}\frac{1}{n}\cdot\frac{\vert2^n-1\vert_3}{2^n}
+\bigo\left(2^{-N/2}\right)\\
&=&\sum_{n=1}^{\infty}\frac{1}{n}\cdot\frac{\vert2^n-1\vert_3}{2^n}+\bigo\left(2^{-N/2}\right).
\end{eqnarray*}
In particular, the difference
between~$F(N)$ and the sum in~\eqref{Mersum} is a constant
plus~$\bigo\left(2^{-N/2}\right)$.

Some well-known partial sums related to the classical Mertens'
Theorem will be needed. For~$x>0$,
\[
\sum_{n \le x}\frac{1}{n}=\log x + c_2 +\bigo(1/x),
\]
where the constant~$c_2$ is the Euler--Mascheroni constant. It
follows that
\begin{equation}\label{uniform}
\sum_{k\le x,\gcd(p,k)=1}\frac{1}{k}=\left(\frac{p-1}{p}\right)\log
x + c_3(p) + \bigo(1/x)
\end{equation}
for any prime~$p$, where~$c_3(p)$ is a constant depending on~$p$
(the implied constant in the~$\bigo(1/x)$ term also depends on~$p$).

The sum in~\eqref{Youweretheonewhobroughtmehomedowntoearth} can be
estimated using~\eqref{KatieMelua} as follows. The sum over the odd
terms is
\[
\sum_{n \le N,2\smallnotdivides n}\frac{1}{n}=\frac{1}{2}\log N +
c_4 + \bigo(1/N)
\]
by~\eqref{uniform}, with~$c_4=c_3(2)$. The sum over the even terms
collapses just as before to give
\[
\sum_{2k \le N}\frac{\vert 3k\vert _3}{2k}.
\]
Now
\[
\displaystyle \sum_{k \le N}\frac{\vert k\vert
_3}{k}
=
\sum_{r=0}^{{\log N}/{\log
3}}\frac{1}{3^{2r}}\sum_{k=1,\gcd(3,k)=1}^{{N}/{3^r}}\frac{1}{k}.
\]
By~\eqref{uniform}, this is
\[
\sum_{r=0}^{{\log N}/{\log 3}}\frac{2}{3^{2r+1}}[\log N - r\log 3 +
c_6 + \bigo(3^r/N)],
\]
where the constant in the~$\bigo(3^r/N)$ term is
independent of~$r$. The
computation of each term involves summing a geometric series. In
each case the sum differs from the full series with an error
that is~$\bigo(1/N)$; we deduce that
\begin{equation*}
\displaystyle \sum_{k \le N}\frac{\vert k\vert
_3}{k}=\frac{3}{4}\log N + c_5 + \bigo(1/N).
\end{equation*}
The sum over the odd and even terms gives
\[
\frac{1}{6}\cdot \frac{3}{4}\log N + \frac{1}{2} \log N + c_7 +
\bigo(1/N)=\frac{5}{8} \log N + c_7 + \bigo(1/N),
\]
completing the proof of Theorem~\ref{theorem:mertensinsimplecases23}
for the case~$\xi=2$ and~$S=\{3\}$.

\section{Proof of
Theorem~\ref{theorem:moregeneralorbitcounting}}\label{section:proofofmoregeneralorbittheorem}

We are given an algebraic number field~$\mathbb K$ with set of
places~$P(\mathbb K)$ and set of infinite places~$P_{\infty}(\mathbb
K)$, an element of infinite multiplicative order~$\xi\in\mathbb
K^*$, and a finite set~$S\subset P(\mathbb K)\setminus
P_{\infty}(\mathbb K)$ with the property that~$\vert\xi\vert_w\le1$
for all~$w\notin S\cup P_{\infty}(\mathbb K)$. The associated ring
of~$S$-integers is
\[
R_S=\{x\in\mathbb K\mid\vert x\vert_w\le1\mbox{ for all }
w\notin S\cup P_{\infty}(\mathbb K)\}.
\]
The compact group~$X$ is the character group of~$R_S$,
and the endomorphism~$T$ is the dual of the map~$x\mapsto\xi x$
on~$R_S$. Examples of this construction may be found in~\cite{MR99b:11089}.
Following Weil~\cite[Chap.~IV]{MR96c:11002}, write~$\mathbb K_w$
for the completion at~$w$, and for~$w$ finite, write~$r_w$ for the
maximal compact subring of~$\mathbb K_w$.

Define the compact group~$X^*_w$ by
\[
X^*_w=
\begin{cases}
\mathbb S^1&\mbox{if }w\in P_{\infty}(\mathbb K)\mbox{ and }
\vert\xi\vert_w=1;\\
r_w^*&\mbox{if }w\notin P_{\infty}(\mathbb K)\mbox{ and }
\vert\xi\vert_w=1;\\
\{1\}&\mbox{in all other cases.}
\end{cases}
\]
Finally, let~$X^*=\prod_{w}X_w^*.$ The element~$a_T=(a_{T,w})_w$
of~$X^*$ is defined by~$a_{T,w}=\imath_w(\xi)$ where~$\imath_w$ is
the corresponding embedding of~$\mathbb K$ into~$\mathbb C$
or~$\mathbb K_w$ whenever~$X_w^*$ is
non-trivial, and~$a_{T,w}=1$ in all other cases.

By~\cite[Lem.~5.2]{MR99b:11089}, the number of points in~$X$
fixed by~$T^n$ is
\begin{equation}\label{equation:periodicpointformulaingeneral}
\fix_{T}(n)=\prod_{w\in S\cup P_{\infty}(\mathbb K)}\vert
\xi^n-1\vert_w,
\end{equation}
so the number of
orbits of
length~$n$ is
\[
\orbit_{T}(n)={\frac{1}{n}}\sum_{d\vert n}\mu(
{\textstyle\frac{n}{d}})
\prod_{w\in S\cup P_{\infty}(\mathbb K)}\vert
\xi^n-1\vert_w
\]
by M{\"o}bius inversion, and hence
\begin{equation}\label{equation:generalcaseexpressionforpi}
\pi_{T}(N)=\sum_{n\le N}{\frac{1}{n}}\sum_{d\vert n}\mu(
\textstyle\frac{n}{d}) \displaystyle\prod_{w}\vert \xi^n-1\vert_w
\end{equation}
where~$w$ is restricted to run through the places in~$S\cup
P_{\infty}(\mathbb K)$ only (both here and below).

We begin by replacing~\eqref{equation:generalcaseexpressionforpi}
with a more manageable expression just as
in~\eqref{equation:GNfirstdefined}. Let
\[
G(N)=\sum_{n\le N}{\frac{1}{n}}
\prod_{\vert\xi\vert_w>1}\vert\xi\vert_w^n
\prod_{\vert\xi\vert_w\le1}\vert\xi^n-1\vert_w.
\]
By~\cite{MR99b:11089}, the topological entropy of~$T$
is
\begin{equation}\label{entropyboundingthingie}
h(T)=\sum_{\vert\xi\vert_w>1}\log\vert\xi\vert_w>0.
\end{equation}
It follows that~$\left(\Pi_T(N)\right)$ is a bounded sequence.
Let~$h'(T)$ denote the maximum value of~$\frac{1}{2}h(T)$ and the
expression~\eqref{entropyboundingthingie} with one term omitted;
notice in particular that~$h'=h'(T)<h=h(T)$. Write
\[
C_{\mathbb K}=4^{\vert P_{\infty}(\mathbb K)\vert}
\]
Now
\begin{eqnarray*}
\left\vert G(N)-\pi_{T}(N)\right\vert&=& \sum_{n\le
N}{\frac{1}{n}}\left( \sum_{d\divides
n}\right.\bigo\left(e^{nh'}\right)
\underbrace{\prod_{\vert\xi\vert_w\le1}\vert\xi^d-1\vert_w}_{\le C_{\mathbb K}}\\
&&
\medspace\medspace\medspace\medspace\medspace\medspace\medspace\medspace
\medspace\medspace\medspace\medspace\medspace\medspace\medspace\medspace
 \left.+ \sum_{d\divides n,d<n} \prod_{w}\vert\xi^n-1\vert_w
\right)\\
&=& \sum_{n\le N}{\frac{1}{n}} \left(\vphantom{\prod}\right.
n\bigo\left(\vphantom{A^A}\right.e^{nh'}\left.\vphantom{A^A}\right)
+\underbrace{\sum_{d\le\lfloor n/2\rfloor}
\prod_{w}\vert\xi^n-1\vert_w}_{\bigo(e^{nh/2})}
\left.\vphantom{\prod}\right)\\
&=& \bigo\left(e^{Nh'}\right).
\end{eqnarray*}
Since~$h'<h$, this means that~$\left(\Pi_T(N_j)\right)$
converges if
and only if
\[
\frac{N_jG(N_j)}{e^{h(T)(N_j+1)}}
\]
converges. Write
\[
G(N)=\sum_{n\le N}{\frac{1}{n}} A(n)B(n)
\]
where
\[
A(n)=\prod_{\vert\xi\vert_w>1}\vert\xi\vert_w^n,
\]
and
\[
B(n)=\prod_{\vert\xi\vert_w\le1}\vert\xi^n-1\vert_w.
\]
Notice that~$A(n)=e^{hn}$,~$B(n)\le C_{\mathbb K}$, and a
subsequence~$\left(B(N_j)\right)$ of~$\left(B(N)\right)$ converges
whenever~$\left(a_T^{N_j}\right)$ converges in~$X^*$ (since the
terms in~$B(N)$ with~$\vert\xi\vert_w<1$ simply converge to~$1$).

As before, let~$K(N)=\lfloor N^{1/4}\rfloor$, and consider the
expression
\begin{eqnarray*}
    a_N&=&\sum_{n=N-K(N)}^{N}\frac{N}{e^{h(N+1)}}\cdot\frac{1}{n}\cdot
    A(n)B(n)\\
&=&\sum_{t=0}^{K(N)}\frac{N}{e^{h(N+1)}}\cdot\frac{1}{N-t}A(N-t)B(N-t).
\end{eqnarray*}
Now
\begin{eqnarray}
  \left\vert a_N-\frac{G(N)N}{e^{h(N+1)}}\right\vert &=&
  \sum_{t=K(N)+1}^{N}\frac{NA(N-t)B(N-t)}{(N-t)e^{h(N+1)}}\nonumber \\
  &\le& \sum_{t=K(N)+1}^{N}\frac{N\cdot C_{\mathbb K}}{e^{h(t+1)}}\nonumber \\
  &=& \bigo\left(Ne^{-K(N)}\right)\label{equation:estimate1}
\end{eqnarray}
so in order to show that~$\left(\Pi_T(N_j)\right)$ converges it is
enough to show that the subsequence~$(a_{N_j})$ converges. The
expression for~$a_N$ can be further simplified, since
\begin{eqnarray}
  a_N &=& \sum_{t=0}^{K(N)}\frac{N}{e^{h(N+1)}}\cdot\frac{1}{N-t}A(N-t)B(N-t)\nonumber \\
   &=& \sum_{t=0}^{K(N)}\frac{1}{e^{h(t+1)}}\cdot\frac{1}{1-t/N}B(N-t)\nonumber \\
   &=& a_N^*+{\bigo\left(\sum_{t=0}^{K(N)}
   \frac{t}{N}C_{\mathbb K}\right)}=a_N^*+{\bigo\left(N^{-1/2}\right)},\label{equation:estimate2}
\end{eqnarray}
where
\[
a_{N}^*=\sum_{t=0}^{K(N)}\frac{1}{e^{h(t+1)}}\cdot B(N-t).
\]
Choose~$\delta$ with
\[
0<\delta={\textstyle\frac{1}{2}}\min\{\vert\xi^j-1\vert_w\mid
\vert\xi\vert_w=1,1\le j\le\vert S\vert, w\in S\}.
\]
If~$\vert\xi^N-1\vert_w<\delta$, then
\[
\vert\xi^{N-j}-1\vert_w=\vert\xi^{-j}(\xi^N-1)+\xi^{-j}-1\vert_w
\ge\vert\xi^{-j}-1\vert_w-\delta>\delta \mbox{
for }1\le j\le\vert S\vert.
\]
%
%
Notice that~$a_N^*$ can only be small if~$B(N),B(N-1),\dots,B(N-\vert S\vert)$
are small, but if
\[
\vert B(N-j)\vert<\delta^{\vert S\vert}\mbox{ for }j=0,\dots,\vert S\vert-1
\]
then~$\vert B(N)-\vert S\vert)\vert>\delta^{\vert S\vert}$.
It follows that there is no sequence~$(N_j)$ with
\[
\prod_{\vert\xi\vert_w\le1}\vert\xi^{N_j+k}-1\vert_w\rightarrow0\mbox{
for }k=0,1,2,\dots,
\]
and, indeed~$\liminf_{N\to\infty}a_{N}^*\ge\delta^{\vert S\vert}>0$.

Assume now that~$(N_j)$ is a sequence with the property
that~$\left(a_T^{N_j}\right)$ converges in~$X^*$, so in particular each
sequence~$\left(\vert\xi^{N_j}-1\vert_w\right)$ is Cauchy for~$w\in S$,~$\vert\xi\vert_w\le1$,
hence~$\left(\vert\xi^{N_j-t}-1\vert_w\right)$
and~$\left(B(N_j-t)\right)$ are Cauchy for each~$t$.
Moreover, these sequences are uniformly Cauchy in~$t$,
since~$\vert\xi^{N_j-t}-\xi^{N_k-t}\vert_w=\vert\xi^{N_j}-\xi^{N_k}\vert_w$
for all~$t$.
We claim
that~$(a_{N_j}^*)$ also converges, which (by the
estimates~\eqref{equation:estimate1} and~\eqref{equation:estimate2})
will complete the proof of
Theorem~\ref{theorem:moregeneralorbitcounting}. Let~$k<j$ be fixed.
Then\begin{eqnarray*}
  \vert a_{N_j}^*-a_{N_k}^*\vert &\le& \left\vert
  \sum_{t=0}^{K(N_j)}\frac{1}{e^{h(t+1)}}B(N_j-t)\right.\\
  &&\medspace\medspace\medspace\medspace\medspace\medspace\medspace\medspace
  \left.- \sum_{t=0}^{K(N_k)}\frac{1}{e^{h(t+1)}}B(N_k-t)
  \right\vert \\
   &\le& \sum_{t=0}^{K(N_k)}\frac{1}{e^{h(t+1)}}{\left\vert
   B(N_j-t)-B(N_k-t)\right\vert}\\
  &&\medspace\medspace\medspace\medspace\medspace\medspace\medspace\medspace
  +\underbrace{\sum_{t=K(N_k)+1}^{K(N_j)}\frac{1}{e^{h(t+1)}}B(N_j-t)}_{\bigo\left(e^{-hK(N_k)}\right)}\\
  &\longrightarrow&0\mbox{ as }k\to\infty,
\end{eqnarray*}
since
\[
\sum_{t=0}^{K(N_k)}\frac{1}{e^{h(t+1)}}\left\vert
   B(N_j-t)-B(N_k-t)\right\vert
   \]
   \[
\le\left(\sum_{t=0}^{\infty}\frac{1}{e^{h(t+1)}}\right)
   \underbrace{\max_{0\le t\le K(N_k)}\left\vert
   B(N_j-t)-B(N_k-t)\right\vert}_{\to0\text{ as }k\to\infty\text{ by the uniform
Cauchy property}}.
\]


\section{Mertens' Theorem without error
term}\label{section:mertenswithouterrorterm}

The setting is an~$S$-integer map~$T:X\to X$ with~$X$ connected
and~$S$ finite. We first give a simple argument to show a form of
Theorem~\ref{theorem:mertensforSfinite} without error term, and then
consider how an error term is obtained. Recall that
\[
    \mertens_T(N)=\sum_{n\le N}\frac{1}{n}\sum_{d\divides
    n}\mu(n/d)\left(\frac{\prod_{w}\vert\xi^d-1\vert_w}{e^{hn}}\right).
\]
Let
\[
C(n)=\prod_{\vert\xi\vert_w\neq 1}\vert\xi^n-1\vert_w
\]
and
\[
D(n)=\prod_{\vert\xi\vert_w=1}\vert\xi^n-1\vert_w.
\]
Define
\[
    F(N)=\sum_{n\le
    N}\frac{1}{n}D(n),
\]
and write
\[
h^*=\prod_{\genfrac{}{}{0pt}{}{\vert\xi\vert_w>1,}{w\vert\infty}}\vert\xi\vert_w
\]
for the Archimedean contribution to the entropy. Then
\begin{eqnarray*}
\mertens_T(N)-F(N)&=&\sum_{n\le N}\frac{1}{n}\left(\sum_{d\divides
n}\mu\left({\textstyle\frac{n}{d}}\right)e^{-hn}\prod_{w}\vert\xi^d-1\vert_w-D(n)\right)\\
&=&\sum_{n\le N}\frac{1}{n}\sum_{d\divides
n}\mu\left({\textstyle\frac{n}{d}}\right)D(d)
\prod_{\genfrac{}{}{0pt}{}{\vert\xi\vert_w>1,}{w\vert\infty}}{\vert\xi^d-1\vert_w}{\vert\xi\vert_w^n}\\
&&\medspace\medspace\medspace\medspace\medspace\medspace\medspace\medspace
\medspace\medspace\medspace\medspace\medspace -\sum_{n\le
N}\frac{1}{n}D(n)\\
&=&\sum_{n\le
N}\frac{1}{n}\left({\vphantom{\sum}}\right.D(n)\left(1-\bigo(e^{-h^*n})\right)-D(n)\left.{\vphantom{\sum}}\right)\\
&&\medspace\medspace\medspace\medspace\medspace\medspace\medspace\medspace
\medspace\medspace\medspace\medspace\medspace+\sum_{n\le
N}\frac{1}{n}\bigo\left({\vphantom{\sum}}\right.\underbrace{\sum_{d<n/2}D(d)e^{h^*(d-n)}}_{\bigo(e^{-h^*n/2})}
\left.{\vphantom{\sum}}\right)\\
&=& \sum_{n\le N}\frac{1}{n}D(n)\bigo(e^{-h^*n})+\sum_{n\le
N}\frac{1}{n}\bigo(e^{-h^*n/2})
\end{eqnarray*}
in which the implied constants are uniformly bounded. It follows
that~$\mertens_T(N)-F(N)$ may be written as the difference between a
sum of a convergent series and the sum from~$N$ to~$\infty$ of that
series, and this tail of the series is~$\bigo(e^{-h^*N})$. Thus in
order to prove Theorem~\ref{theorem:mertensforSfinite} it is enough
to consider~$F(N)$.

\begin{lemma}\label{kroneckerweil}
Let~$g$ be an element of a compact abelian group~$G$.
Then the sequence~$\left(g^n\right)$ is uniformly
distributed in the smallest closed subgroup of~$G$
containing~$g$.
\end{lemma}

\begin{proof}
This is essentially the Kronecker--Weyl lemma. Write~$X$ for the
closure of the set~$\{g^n\mid n\in\mathbb Z\}$ and~$\mu_X$ for Haar
measure on~$X$. In order to show that
\[
\frac{1}{N}\sum_{n=1}^{N}f(g^n)\rightarrow\int f\dee\mu_X
\]
for all continuous functions~$f:X\to\mathbb C$, it is
enough to show this for characters. If~$\chi:X\to\mathbb S^1$
is a non-trivial character on~$X$, then
\begin{eqnarray*}
\left\vert\frac{1}{N}\sum_{n=0}^{N-1}\chi(g^n)\right\vert
&=&
\left\vert\frac{1}{N}\sum_{n=0}^{N-1}\chi(g)^n\right\vert\\
&=&
\left\vert\frac{1}{N}\cdot\frac{1-\chi(g)^N}{1-\chi(g)}\right\vert\\
&\le&\frac{1}{N}\cdot\frac{2}{1-\chi(g)}\rightarrow0\mbox{ as }N\to\infty,
\end{eqnarray*}
so the sequence is uniformly distributed.
\end{proof}

Lemma~\ref{kroneckerweil} may be applied to the element~$a_T\in
X^*$: the function
\[
x\mapsto\prod_{\vert\xi\vert_w=1}\vert x-1\vert_w
\]
is continuous on~$X^*$, so
\[
\frac{1}{N}\sum_{n=1}^{N}D(n)\rightarrow
k_T\mbox{ as }N\to\infty
\]
where
\[
k_T=\int_{X^*}\prod_{\vert\xi\vert_w=1} \vert
x-1\vert_w\dee\mu_{X^*}.
\]
Thus
\begin{eqnarray*}
F(N)&=&\sum_{n=1}^{N}\left(
\frac{1}{n}-\frac{1}{n+1}\right)\sum_{m=1}^{n}D(m)+\frac{1}{N+1}
\sum_{m=1}^{N}D(m)\\
&\sim&k_T\log N,
\end{eqnarray*}
giving Theorem~\ref{theorem:mertensforSfinite} without error term.

\section{Mertens Theorem with~$\mathbb K=\mathbb
Q$}\label{section:mertenswithkequalsq}

Section~\ref{section:boththeoremsfor23example} contains a proof of
Theorem~\ref{theorem:mertensinsimplecases23} for
the case~$S=\{3\}$ and~$\xi=2$. In this section we prove
Theorem~\ref{theorem:mertensinsimplecases23}; the essential
difference between this and Theorem~\ref{theorem:mertensforSfinite}
is that the assumption~$\mathbb K=\mathbb Q$ does not permit~$\xi$
to induce an ergodic map (that is,~$\xi$ is not a
unit root) while exhibiting non-hyperbolicity in an
infinite place. The argument in this section, with simple
modifications, would give
Theorem~\ref{theorem:mertensinsimplecases23} under the assumption
that~$\mathbb K$ does not contain any Salem numbers ($[\mathbb
K:\mathbb Q]\le3$ would suffice, for example).

Fix a finite set~$S$ of primes, a rational~$r\in\mathbb Q$
with~$r\neq\pm1$ and~$\vert r\vert_p<1\implies p\in S$. Consider the
map~$T:X\to X$ dual to the map~$x\mapsto rx$ on the additive group
of the ring
\[
R_S=\{r\in\mathbb Q\mid\vert r\vert_p\le1\mbox{ for all }p\notin
S\}.
\]
By~\cite[Lem.~5.2]{MR99b:11089}, the number of points fixed by~$T^n$
is
\[
\fix_{T}(n)=\vert r^n-1\vert\prod_{p\in S}\vert
r^n-1\vert_p=(r^n-1)\vert r^n-1\vert_S,
\]
where we write~$\vert x\vert_S$ for~$\prod_{p\in S}\vert x\vert_p$,
and so
\[
\orbit_{T}(n)={\frac{1}{n}}\sum_{d\vert n}\mu(
{\textstyle\frac{n}{d}})\vert r^d-1\vert|r^d-1|_S.
\]
Just as in Section~\ref{section:mertenswithouterrorterm},
it is sufficient to work with the sum~$F(N)$.

The analog of Mertens' Theorem in this setting is most easily proved
by isolating the following arithmetic argument. A function~$f$ is
called totally multiplicative if~$f(mn)=f(m)f(n)$ for
all~$m,n\in\mathbb N$.

\begin{lemma}\label{inclusionexclusionlemma}
Let~$f:\mathbb N\to\mathbb C$ be a totally multiplicative function
with
\[
\sum_{n\le N} f(n) = k_f\log N + c_f + \bigo(1/N),
\]
for constants~$c_f$ and~$k_f$. Let~$E$ be a finite set of natural
numbers and, for~$D\subseteq E$, let~$n_{D}=
\lcm\{n:n\in D\}$. Then there is a
constant~$c_{f,E}$ for which
\[
\sum_{n\le N,k\smallnotdivides n\,{\rm for}\,k\in E} f(n) =
k_{f,E}\log N+c_{f,E} +\bigo(1/N),
\]
where
\[
k_{f,E}=k_f \sum_{D\subseteq E} (-1)^{\vert D\vert} f(n_{E}).
\]
\end{lemma}

\begin{proof}
Notice that
\begin{eqnarray*}
\sum_{n\le N, n_{D}\divides n} f(n) &=&
f(n_{D}) \sum_{n\le N/n_{D}} f(n) \\
&=& f(n_{D})
\left( k_f \log(N/n_{D}) + c_f + \bigo(1/N) \right) \\
&=& k_f f(n_{D})\log N + c_{f,n_D} + \bigo(1/N),
\end{eqnarray*}
for some constant~$c_{f,n_{D}}$. The result follows by an
inclusion-exclusion argument.
\end{proof}

Notice that, if~$E$ is a set of pairwise coprime natural numbers,
then
\[
k_{f,E} = k_f \prod_{n\in E} \left(1-f(n)\right).
\]

Now let~$\cP$ be a finite set of (rational) primes.
For~$\br=(r_p)_{p\in\cP}\in\mathbb Z^{|\cP|}$, write
\[
\bp^\br = \prod_{p\in\cP} p^{r_p}
\]
and abbreviate~$\bp=\bp^{(1,\dots,1)}=\prod_{p\in\cP}p$. Define a
partial order on~$|\cP|$-tuples by
\[
\br=(r_p)_{p\in\cP}\le\bs=(s_p)_{p\in\cP} \ \iff\ r_p\le s_p\
\forall p\in\cP
\]
and write~$\bO=(0)_{p\in\cP}$.

For~$\bt=(t_p)_{p\in\cP}\in\mathbb N^{|\cP|}$, write
\[
f_{\cP,\bt}(n)=\frac 1n \prod_{p\in\cP}|n|_p^{t_p};
\]
notice that this is a totally multiplicative function.

\begin{proposition} \label{Mertprop}
There is a constant~$c_{\cP,\bt}$ for which
\[
\sum_{n<N} f_{\cP,\bt}(n) = k_{\cP,\bt} \log N +c_{\cP,\bt}+
\bigo(1/N).
\]
where~$k_{\cP,\bt}$ is the
product~$\prod_{p\in\cP}\left(1-\frac 1p\right)\left(1-\frac
1{p^{t_p+1}}\right)^{-1}$.
\end{proposition}

Note that, since~$f_{\cP,\bt}$ is totally multiplicative,
Lemma~\ref{inclusionexclusionlemma} may be applied to this result to get asymptotics
for sums over subsets of~$\mathbb N$.

\begin{proof} The proof is by induction on~$m=|\cP|$, the case~$m=0$
being the familiar statement
\[
\sum_{n<N}\frac{1}n = \log N +c + \bigo(N^{-1}).
\]
Write~$p^r\Divides n$ if~$r=\ord_p(n)$ is the exact order with
which~$p$ divides~$n$.
Put~$\cP=\{p_1,...,p_m\}$,~$\cP_1=\{p_2,...,p_m\}$,~$t_1=t_{p_1}$
and~$\bt_1=(t_{p_2},...,t_{p_m})$. Then
\begin{eqnarray*}
\qquad\sum_{n\le N} f_{\cP,\bt}(n)&=& \sum_{r_1=0}^{\log N/\log p_1}
\sum_{n\le N, p_1^{r_1}\Divides n} f_{\cP,\bt}(n)
\\
&=& \sum_{r_1=0}^{\log N/\log p_1} {\frac 1{p_1^{(t_1+1)r_1}}}
\sum_{n\le N/p_1^{r_1},p_1\smallnotdivides n} f_{\cP_1,\bt_1}(n) \\
&=& \sum_{r_1=0}^{\log N/\log p_1} \frac 1{p_1^{(t_1+1)r_1}}
\left(1-{\frac 1{p_1}}\right)k_{\cP_1,\bt_1} \cdot \\
&&\qquad\qquad\qquad \left[\log N -  r_1\log p_1 + c' +
\bigo(p_1^{r_1}/N)\right]
\end{eqnarray*}
using the inductive hypothesis and
Lemma~\ref{inclusionexclusionlemma} (applied to~$f=f_{\cP_1,\bt_1}$
and~$E=\{p_1\}$). Note that the implied constants in
the~$\bigo(p_1^{r_1}/N)$ terms are
independent of~$r_1$. The computation
of each term involves summing some geometric series, and in each
case the sum differs from the full series with an error term
that is~$\bigo(1/N)$.
\end{proof}

The next argument will be needed again in
Section~\ref{section:allowinginfiniteplaces} in a more general setting,
so we now allow~$\mathbb K$ to be a number field.
Theorem~\ref{theorem:mertensinsimplecases23} will follow
at once, since the sum considered here is the~$F(N)$ from
Section~\ref{section:mertenswithouterrorterm}.

\begin{proposition}\label{NearlyMert}
Let~$\mathbb K$ be a number field,~$\xi\in\mathbb K$ and~$S$ a
finite set of
non-Archimedean places of~$\mathbb K$ such that~$|\xi|_\fp=1$ for
all~$\fp\in S$. Write~$|x|_S=\prod_{\fp\in S}|x|_\fp$
for~$x\in\mathbb K$. Then there are constants~$k_S\in\mathbb Q$
and~$c_S\in\mathbb R$ such that
\[
\sum_{n<N}\frac{|\xi^n-1|_S}n = k_S\log N +c_S+ \bigo(1/N).
\]
\end{proposition}

\begin{proof}
For~$\fp\in S$, let~$o_\fp$ denote the order of~$\xi$ in the residue
field at~$\fp$, that is, the least positive integer~$o$ such
that~$|\xi^o-1|_\fp<1$. Then
\[
|\xi^n-1|_\fp=1 \ \iff\ o_\fp\notdivides n.
\]
Let~$p$ be the rational prime such that~$\fp\divides p$. It is
sometimes more convenient to use the extension of the~$p$-adic
absolute value~$|\cdot|_p$, which is related to~$|\cdot|_\fp$ by
\[
|\cdot|_\fp^{\vphantom{[K_\fp:\mathbb Q_p]}} =
|\cdot|_p^{[K_\fp:\mathbb Q_p]},
\]
where~$K_\fp$ is the completion of~$K$ at~$\fp$.

Let~$m_\fp$ be the least positive integer~$m$ such that
\[
|\xi^m-1|_p < \frac 1{p^{1/p-1}}.
\]
Then~$m_\fp=p^{r_\fp}o_\fp$, for some~$r_\fp\ge 0$. Moreover,
if~$m_\fp|n$ then
\[
|\xi^n-1|_\fp = |n|_\fp |\log\xi|_\fp,
\]
where~$\log$ is here the~$p$-adic logarithm.

Finally, if~$n=k p^r o_\fp$, with~$(k,p)=1$, then
\[
|\xi^n-1|_\fp = |\xi^{p^r o_\fp}-1|_\fp.
\]

For~$T$ a subset of~$S$, put~$o_T=\lcm\{o_\fp:\fp\in T\}$. Split up
the sum according to the subsets of~$S$, giving
\[
\sum_{n<N}\frac{|\xi^n-1|_S}n = \sum_{T\subset S}\ \sum_{n<N,\,
o_T|n,\, o_\fp\smallnotdivides n\forall \fp\not\in T}
\frac{|\xi^n-1|_T}n.
\]
We show that each internal sum has the required form and, since
there are only a finite number of subsets of~$S$, we will be done.

So let~$T$ be a subset of~$S$ and let~$\cP$ be the set of rational
primes divisible by some~$\fp\in T$. Putting~$m_T=\lcm\{m_\fp:\fp\in
T\}$, there exists~$\br=(r_p)\ge\bO$ such that~$m_T=\bp^\br o_T$.
Then we have
\begin{eqnarray*}
&\hskip-20pt& \sum_{n<N,\, o_T|n,\, o_\fp\smallnotdivides n\forall
\fp\not\in T} \frac{|\xi^n-1|_T}n
\\
&\hskip-20pt&\hskip70pt=\ 
\sum_{\bO\le\bs\le\br} \frac{|\xi^{\bp^\bs o_T}-1|_T}{\bp^\bs o_T}
\sum_{n<N/\bp^\bs o_T,\, (n,\bp)=1,\, o_\fp\smallnotdivides n\bp^\bs
o_T\forall \fp\not\in T}
\frac{1}n \\ \\
&\hskip-20pt&\hskip120pt+\ \frac{|\xi^{m_T}-1|_T}{m_T}
\sum_{n<N/m_T,\,o_\fp\smallnotdivides nm_T\forall \fp\not\in T}
\frac{|n|_T}n.
\end{eqnarray*}
Now~$|n|_T=\prod_{p\in\cP}|n|_p^{t_p}$, where~$t_p=\sum_{\fp\in
T,\,\fp|p}[K_\fp:\mathbb Q_p]$, so~$\frac{|n|_T}n=f_{\cP,\bt}(n)$.
So this again gives a finite number of sums, each of which has the
required form, by applying Lemma~\ref{inclusionexclusionlemma} to
Proposition \ref{Mertprop}.
\end{proof}

This completes the proof of
Theorem~\ref{theorem:mertensinsimplecases23}. The constants
appearing in Theorem~\ref{theorem:mertensinsimplecases23} may be
found explicitly for any given set~$S$, by following
the recipe in the proof of Proposition~\ref{NearlyMert} and using
Proposition~\ref{Mertprop}, leading to
Example~\ref{example:explicitconstants}.

\section{Allowing infinite
places}\label{section:allowinginfiniteplaces}

The estimate in~\eqref{equation:threetermestimate} requires several
improvements to the argument above. From now on~$S$ denotes a finite
set of non-Archimedean valuations on the number field~$\mathbb K$
and~$\xi\in\mathbb K^*$ is an element of infinite multiplicative
order with~$\vert\xi\vert_v=1$ for all~$v\in S$.

\begin{lemma}\label{equalM}
Let~$M\in \mathbb N$ denote any integral~$S$-unit. The solutions of
the equation
\[
|\xi ^n-1|_S=\frac{1}{M}
\]
consist of~$\bigo(M^{1-1/d})$ cosets mod~$M'$ where~$M'=\rho M$ for
some fixed integer~$\rho$ and some~$d>0$, both independent of~$M$.
\end{lemma}

\begin{proof}
For each~$v\in S$, the set~$U_k=\{n\in\mathbb Z\mid\ord_v(\xi^n-1)\ge k\}$ is
a subgroup of~$\mathbb Z$. For sufficiently large~$k$, the cosets
of~$U_{k+1}$ in~$U_k$ are defined by either~$1$ or~$p$ congruence
classes modulo~$sp^{k+1}$ for a uniform constant~$s$. Now for~$n\in
U_{k}\setminus U_{k+1}$,~$\vert n\vert_v=sp^{-kd}$ for~$d=[\mathbb K:\mathbb Q]$,
so~$n$ lies in~$\bigo(p^{kd-k})=\bigo(M^{1-1/d})$ classes.
Choose~$\rho=m_v$ in the notation of the proof of Proposition~\ref{NearlyMert}.
The
Chinese Remainder Theorem then gives the same
bound for the product of the finitely many valuations in~$S$.
\end{proof}

Write~${\sum}'$ for a sum taken only over
integral~$S$-units.

\begin{lemma}
For any~$c>0$, the series
\begin{equation}\label{sumoverM}
{\sum_M}'\frac{\log M}{M^c}
\end{equation}
converges. The tail of the series satisfies
\[
{\sum_{M>X}}'\frac{\log M}{M^c}=\bigo(1/X^e),
\]
for any~$e<c$.
\end{lemma}

\begin{proof}
Let~$p_1,\dots,p_r$ be the distinct rational primes dividing
the elements of~$S$.
Write each integral~$S$-unit~$M$ in the form~$p_1^{e_1}\dots p_r^{e_r}$
with~$0\le e_i$ for~$i=1,\dots ,r$. The sum
in~\eqref{sumoverM} is then a finite sum of terms, each of which may
be written as a finite product of convergent geometric progressions
and their squares, showing the convergence. To estimate the error
notice that if~$M>X$ then at least one term~$e_i>\kappa\log X$ for
some uniform constant~$\kappa$, depending on~$S$ only. Hence the
error is bounded above by
\[
\sum_{i=1}^rK_i\sum_{t>\kappa\log X}\frac{t}{p_i^{ct}},
\]
for some constants~$K_i$, and this sum is~$\bigo(\log X/X^c)$ by Euler
Summation.
\end{proof}

\begin{theorem}\label{finaltheorem?}
Let~$a$ denote a complex algebraic number with~$|a|=1$
and~$a$ not a root of unity. Then for some~$\delta>0$ and constant~$\ell$,
\[
\sum_{n<N}\frac{a^n|\xi ^n-1|_S}{n}=\ell+\bigo(N^{-\delta}).
\]
\end{theorem}

\begin{proof}
Decompose the sum according to the integral~$S$-units~$M$
with
\[
\vert\xi^n-1\vert_S=\frac{1}{M}.
\]
Consider the sum
\[
F_N(X)={\sum_{M<X}}'\frac{1}{M}\sum_{n<N:|\xi
^n-1|_S=\frac{1}{M}}\frac{a^n}{n}.
\]
We claim that there is a constant~$\ell$ for which
\begin{equation}\label{claim}
F_N(X)=\ell+\bigo(\max \{X^B/N,1/X^e\}),
\end{equation}
where~$e>0$ is a constant depending on~$S$ and~$\xi$ only and~$B$ is
a constant depending on~$\xi$ only. To see this, we use
Lemma~\ref{equalM}: Let~$\{\alpha_i\}$ be
representatives for the~$\bigo\left(M^{1-1/d}\right)$ cosets modulo~$M'=
\rho M$ which are solutions to~$\vert\xi^n-1\vert_S=\frac{1}{M}$.
Then each of the sums
\[
\sum_{n<N:n\equiv \alpha_i\pmod{M'}}\frac{a^n}{n}
\]
can be written using Dirichlet characters in the form
\[
\sum_{n<N}\sum_{j=1}^{M'}c_{ij}\frac{\zeta_j^na^n}{n}
\]
where~$|c_{ij}|= 1/M'$
and each~$\zeta_j$ is an~$M'$th root of
unity (see Apostol~\cite[Chap.~6]{MR0434929} for example).
We can rearrange this double sum to get
\[
\sum_{j=1}^{M'}c_{ij}\sum_{n<N}\frac{\zeta_j^na^n}{n}.
\]
The inner sum is a partial sum of a convergent power series for the
logarithm since~$\zeta_ja\neq1$ (convergence to the logarithm is an
instance of Abel's Theorem; see~\cite[Th.~2.6.4]{cca}). Thus
\[
\sum_{n<N:n\equiv \alpha_i\pmod{M'}}
\frac{a^n}{n}=-\sum_{j=1}^{M'}c_{ij}\log(1-\zeta_ja)+
\sum_{j=1}^{M'}c_{ij}\sum_{n>N}\frac{\zeta_j^na^n}{n}.
\]
Applying Abel Summation to the last sum
gives
\[
\sum_{n<N:n\equiv \alpha_i\pmod{M'}}
\frac{a^n}{n}=-\sum_{j=1}^{M'}c_{ij}\log(1-\zeta_ja)+\bigo\left(\frac{1}{N\min_{j}|1-\zeta_ja|}\right),
\]
using the bound~$|c_{ij}|\le 1/M'$. Thus
the sum sought is
\begin{eqnarray}
\nonumber F_N(X)&=&-\sum_{M<X}\frac{1}{M}\sum_{\alpha
_i}\sum_{j=1}^{M'}c_{ij}\log(1-\zeta_ja)\\
&&\medspace\medspace\medspace\medspace\medspace\medspace
+\sum_{M<X}\sum_{\alpha_i}\frac{1}{M}\bigo\left(\frac{1}{N\min_{j}|1-\zeta_ja|}\right)\label{bothbitsneedbaker}
\end{eqnarray}
in which there are~$\bigo(M^{1-1/d})$ terms~$\alpha_i$.

Both sums
in~\eqref{bothbitsneedbaker} require a lower bound for~$\vert
1-\zeta a\vert$ for~$\zeta$ an~$M'$th root of unity. A bound of
the form~$|1-\zeta a|>A/M'^{B}$ for constants~$A,B>0$ when~$\zeta$
is an~$M'$th root of unity follows from Baker's Theorem~\cite{MR0422171}:
writing~$a=e^{2\pi i\theta}$ and~$\zeta=e^{2\pi ij/M'}$,
the quantity~$\vert
1-e^{2\pi ij/M'}e^{2\pi i\theta}\vert$ is small if and only
if~$\frac{j}{M'}+\theta$
is close to some integer~$K$, in which case~$e^{2\pi i(j/M'+\theta)}-1$
is close to~$2\pi i\left(\frac{j}{M'}+\theta-K\right)$;
by Baker's Theorem there are constants~$A,C>0$ with

\[
\vert M'\log(e^{2\pi ij/M'})-M'\log e^{2\pi i\theta}\vert=\vert 2\pi
iR-M'\log a\vert\ge\frac{A}{M'^C}
\]
for any choice of branches of the logarithm (here~$R-j\in M'\mathbb Z$). It follows that
there are constants~$A,B>0$ with~$|1-\zeta a|>A/M'^{B}$.

The first
sum in~\eqref{bothbitsneedbaker} is bounded in
absolute value by
\[
\sum_{M<X}\frac{1}{M}\sum_{\alpha
_i}\sum_{j=1}^{M'}|c_{ij}||\log(1-\zeta_ja)|
\]
\[
=\bigo\left( \sum_{M<X}\frac{1}{M^{1/d}}\max_{j=1\dots
M'}|\log(1-\zeta_ja)|\right),
\]
using the existence of an absolute bound on the number of
the~$\alpha_i$ from Lemma~\ref{equalM} as  well as the
bound~$|c_{ij}|\le 1/M'$. Thus this term is~$\bigo(\sum_{M<X}\log
M'/M^{1/d})$ and we obtain convergence by comparison with the series
\[
{\sum_{M}}'\frac{\log M}{M^{1/d}}
\]
since~$M'$ and~$M$ are commensurate. Thus at this point, in relation
to~\eqref{claim}, any~$e<1/d$ will do.

To estimate the second sum in~\eqref{bothbitsneedbaker} use Baker's
Theorem in the same way to get an estimate
\[
\bigo\left(\sum_{\alpha_i}\sum_{M<X}\frac{1}{M}.\frac{M'^B}{N}\right)=\bigo(X^B/N).
\]
This concludes the proof of claim~\eqref{claim}. To complete the
proof of Theorem~\ref{finaltheorem?}, note that the sum over those~$n$ with
\[
|\xi ^n-1|_S\le \frac{1}{N^{\epsilon}}
\]
is~$\bigo(N^{-\delta})$ since
\begin{equation*}
\sum_{\vert\xi^n-1\vert_S\le
N^{-\epsilon}}\left\vert\frac{a^n\vert\xi^n-1\vert_S}{n}\right\vert\\
\le N^{-\epsilon}\sum_{n<N}\frac{1}{n}=\bigo(N^{-\delta})\mbox{ for
any }\delta<\epsilon.
\end{equation*}
Thus in estimating the error term, we are allowed to assume that
\[
\frac{1}{M}=|\xi ^n-1|_S >
\frac{1}{N^{\epsilon}}.
\]
In other words, we may write~$X=N^{\epsilon}$ in
claim~\eqref{claim}, where~$\epsilon=\frac{1}{B+1/d}$. This finally
gives an error term~$\bigo(1/N^{\epsilon/d})=\bigo(1/N^{1/dB+1})$.
\end{proof}

As we saw in Proposition~\ref{NearlyMert}, a similar result holds
for the case~$a=1$. We have assembled the material needed to prove
Theorem~\ref{theorem:mertensforSfinite}. By the arguments of
Section~\ref{section:mertenswithouterrorterm} above, it
it enough to show that
\[
F(N)=k_T\log N+C_T+\bigo\left(N^{-\delta}\right)
\]
for some~$\delta>0$, where~$F(N)=\sum_{n<N}\frac{1}{n}D(n)$
and
\begin{eqnarray*}
D(n)&=&\prod_{\vert\xi\vert_w=1}\vert\xi^n-1\vert_w\\
&=&\prod_{\vert\xi\vert_w=1,w\vert\infty}\vert\xi^n-1\vert_w\times
\prod_{\vert\xi\vert_w=1,w<\infty}\vert\xi^n-1\vert_w\\
&=&f(a_1^n,\dots,a_r^n)\times\prod_{\vert\xi\vert_w=1,w<\infty}\vert\xi^n-1\vert_w
\end{eqnarray*}
where~$f$ is an integral polynomial in~$r$ variables,
and~$a_i\in\mathbb S^1$ for~$i=1,\dots,r$ are multiplicatively independent.

This reduces the problem to expressions of the form
\[
\sum_{n<N}\frac{1}{n}a^n\vert\xi^n-1\vert_S
\]
with~$a$ an algebraic number of modulus one that is not a root of unity, to which
Theorem~\ref{finaltheorem?} can be applied, or of the same form
with~$a=1$, to which Proposition~\ref{NearlyMert} may be applied.
Notice in particular that the coefficient of the leading term comes
entirely from the case~$a=1$ covered by
Proposition~\ref{NearlyMert}, and is therefore rational.

\begin{remark}
The leading coefficient in Theorem~\ref{theorem:mertensforSfinite}
can also be described
as~$\lim_{N\to\infty}\frac{1}{N}\sum_{n<N}\vert\xi^n-1\vert_S$,
which is redolent of an integral. There is a sophisticated theory
showing that many~$p$-adic integrals must be rational (see
Denef~\cite{MR902824} for example); is it possible to identify the
limit with an~$S$-adic integral, and is it possible to extend that
theory to handle finitely many valuations?
\end{remark}


\end{document}